\documentclass{IEEEoj}

\usepackage{amsmath,amssymb,amsfonts}
\usepackage{algorithmic}
\usepackage{graphicx,color}
\usepackage{textcomp}

\usepackage{amsthm}
\usepackage[slantedGreek]{newtxmath}
\usepackage[OMLmathsfit]{isomath}
\DeclareMathAlphabet{\mathbfsf}{\encodingdefault}{\sfdefault}{bx}{n}
\usepackage{bm}
\usepackage{envmath}
\usepackage{mathtools}
\usepackage{commath}
\usepackage{siunitx}

\usepackage[caption=false,font=footnotesize]{subfig}

\usepackage{booktabs}
\usepackage{footmisc}  

\usepackage{url}

\newcommand{\mr}{\mathrm}

\newcommand{\veg}[1]{\bm{#1}}     
\newcommand{\mat}[1]{\mathsfbfit{#1}} 
\renewcommand{\vec}[1]{\mathsfbfit{#1}} 
\newcommand{\op}[1]{\mathcal{#1}} 

\newcommand{\n}{\hat{\bm{n}}}

\newcommand\restr[2]{{
        \left.\kern-\nulldelimiterspace 
        #1 
        \vphantom{|} 
        \right|_{#2} 
}}

\newcommand\rst[3]{{
        \left.\kern-\nulldelimiterspace 
        #1 
        \vphantom{|} 
        \right|_{#2}^{#3} 
}}

\usepackage[style=ieee, backend=biber, isbn=false]{biblatex}

\usepackage[backend=biber, style=ieee]{biblatex}
\def\OJlogo{\vspace{-10pt}\includegraphics[height=20pt]{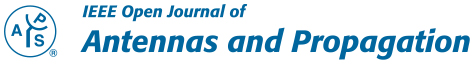}}
\begin{document}
\receiveddate{XX Month, XXXX}
\reviseddate{XX Month, XXXX}
\accepteddate{XX Month, XXXX}
\publisheddate{XX Month, XXXX}
\currentdate{XX Month, XXXX}
\doiinfo{XXXXXXXXXXXX}

\title{Calderón Strategies for the Convolution Quadrature Time Domain Electric Field Integral Equation}

\author{PIERRICK CORDEL\authorrefmark{1}, Student Member, IEEE, ALEXANDRE DÉLY\authorrefmark{1}, ADRIEN MERLINI\authorrefmark{2}, Senior Member, IEEE AND FRANCESCO P. ANDRIULLI\authorrefmark{1},Fellow, IEEE}
\affil{Department of Electronics and Telecommunications (DET), Politecnico di Torino, 10129 Turin, Italy}
\affil{Microwave Department, IMT Atlantique, 29285 Brest, France}
\corresp{CORRESPONDING AUTHOR: F. P. ANDRIULLI (e-mail: francesco.andriulli@polito.it).}
\authornote{This work was supported in part by the European Research Council (ERC) through the European Union’s Horizon 2020 Research and Innovation Programme under Grant 724846 (Project 321) and in part by the H2020-MSCA-ITN-EID project COMPETE GA No 955476.}
\markboth{Preparation of Papers for \textsc{IEEE Open Journal of Antennas and Propagation}}{Author \textit{et al.}}

\begin{abstract}
In this work, we introduce new integral formulations based on the convolution quadrature method for the time-domain modeling of perfectly electrically conducting scatterers that overcome some of the most critical issues of the standard schemes based on the electric field integral equation (EFIE). The standard time-domain EFIE-based approaches typically yield matrices that become increasingly ill-conditioned as the time-step or the mesh discretization density increase and suffer from the well-known DC instability.
This work presents solutions to these issues that are based both on new Calderón strategies and quasi-Helmholtz projectors regularizations. In addition, to ensure an efficient computation of the marching-on-in-time, the proposed schemes leverage properties of the Z-transform---involved in the convolution quadrature discretization scheme---when computing the stabilized operators. The two resulting formulations compare favorably with standard, well-established schemes. The properties and practical relevance of these new formulations will be showcased through relevant numerical examples that include canonical geometries and more complex structures.
\end{abstract}



\maketitle

\section{INTRODUCTION}
\IEEEPARstart{T}{ime} domain boundary integral equations (TDIEs) are widely used in the simulation of transient electromagnetic fields scattered by perfectly electrically conducting (PEC) objects \cite{peterson1998computational,jin2011theory,gibson2021method,colton2013integral}. Like their frequency-domain counterparts, the spatial discretization of these equations is often performed via the boundary element method. The time discretization, however, can be tackled in different ways. A popular approach leverages time basis functions either within a Marching-On-in-Time (MOT) scheme \cite{bennett1968technique,rao1991transient,sayas2016retarded} or within a Marching-On-In-Order procedure \cite{chung2004solution}.
The convolution quadrature (CQ) approach \cite{lubich1988convolutionI,lubich1988convolutionII} is an attractive alternative to these methods in which only space basis functions are explicitly defined. The approach has been applied to several equations in elastodynamics and acoustics \cite{schanz1997new,banjai2012runge} and then in electromagnetics \cite{wang2011implicit}. It provides an efficient time-stepping scheme with matrices derived from the space-discretized Laplace domain operators.

Another advantage of the CQ method is the use of implicit schemes (e.g. Runge Kutta methods \cite{lubich1993runge,oppenheim2001discrete, chipman1971stable,runge1895numerische,kutta1901beitrag}), which are generally more stable and typically allow for a better accuracy control of the solution over time \cite{butcher2016numerical,skvortsov2003accuracy}. However, the CQ time stepping scheme is solved via a computationally expensive MOT algorithm. Nowadays, fast solvers can reach quasi-linear complexity in time and space \cite{banjai2014fast,maruyama2016transient}. Usually, this fast technology uses iterative solvers, resulting in an overall computational cost that is proportional to the number of iterations which is low for well-conditioned systems. Working with well-conditioned matrices is therefore essential to reduce the computational cost of the solution process, in addition to being necessary to obtain accurate results \cite{golub2013matrix}.

Lamentably, however, the CQ discretized time domain electric field integral equation (EFIE) is plagued by several drawbacks. Indeed, the matrices resulting from the discretization of the EFIE are known to become ill-conditioned for large time steps or at dense mesh discretizations: the condition number of the MOT matrices grows quadratically with the time step and with the inverse of the average mesh edge length. These two phenomena are the CQ counterparts of what for standard MOT schemes are known as the large time step breakdown \cite{andriulli2009analysis,chen2001integral,andriulli2012well,bogaert2013low} and the dense discretization breakdown (or $h$-refinement breakdown) \cite{andriulli2009time,cools2009time}. Another challenge in handling the CQ EFIE is that it involves operators whose definitions include a time integration. To avoid dealing with this integral, the time-differentiated counterpart of this formulation is often used \cite{wang2011implicit,dely2019large}, but this differentiation is subject to a source of instability in the form of spurious linear currents living in the nullspace of the operator that degrades the solution \cite{andriulli2009time,shi2014static}; this phenomenon is known as the direct current instability (DC instability). 

In this work, we propose new Calderón-preconditioned and quasi-Helmholtz regularized formulations free from the limitations mentioned above. The Calderón identities they rely on are already a well-established preconditioning approach in both the frequency domain \cite{andriulli2008multiplicative} and time domain discretized by the Galerkin method \cite{andriulli2009time,cools2009time,cools2009nullspaces,beghein2015dc} that is extended in this work to convolution quadrature discretizations and complemented with quasi-Helmholtz regularization. The contribution of this paper is twofold: (i) we present a first approach to tackle the regularization of the EFIE operator and to address the DC instability resulting in a new operator that presents no nullspace on simply connected geometries, thus stabilizing the solution, and (ii) we build upon this first regularized form of the EFIE to obtain an equation that, at the price of a higher number of matrix-vector products, is stable in the case of multiply connected geometries.

This article is structured as follows: the time domain formulations of interest are summarized in Section~\ref{section:Back_Notation} along with the convolution quadrature method and the boundary element method for spatial discretization; in Section~\ref{section:CP_EFIE}, the new Calderón and projectors-based preconditioning strategies are presented; finally, Section~\ref{section:results} presents the numerical studies that confirm the effectiveness of the different approaches before concluding. Preliminary studies pertaining to this work were presented in the conference contribution \cite{cordel2022calderon}.

\section{Background and Notations}
\label{section:Back_Notation}
\subsection{Time Domain Integral Formulations}

In this work, we consider the problem of time-domain scattering by a perfectly electrically conducting object that resides in free space. The object is illuminated by an electromagnetic field $(\veg e^\mr{inc}, \veg h^\mr{inc}) (\veg r,t)$ which induces a surface current density $\veg j_\Gamma$ on its boundary $\Gamma$ that is the solution of the time-domain EFIE
\begin{equation}
        \eta_{0}\op T \left(\veg j_\Gamma \right) \left(\veg r,t\right)= - \n \left(\veg r \right) \times \veg e^\mr{inc} \left(\veg r,t\right)\, .
 \label{eq:EFIE}
\end{equation}
Here, $\n$ is the outpointing normal to $\Gamma$ and $\eta_{0}$ is the characteristic impedance of the background. The electric field operator $\op T$ includes the contributions of the vector and scalar potentials, respectively denoted $\op T_\text{s}$ and $\op T_\text{h}$ \cite{dely2019large}
\begin{gather}
     \op T\left(\veg f \right) \left(\veg r,t \right)=-\frac{1}{c_0}\frac{\partial}{\partial t}\op T_\text{s} \left( \veg f \right)\left(\veg r,t \right) + c_0\int_{-\infty}^{t} \hspace{-7px} \op T_\text{h}\left(\veg f\right)\left(\veg r,t'\right)dt'\,,
\\
  \op T_\text{s}\left(\veg f\right)\left(\veg r,t\right)=\n \left(\veg r \right) \times \iint_{\veg r' \in \Gamma} \left(\op G_{\veg r}*_t\veg f\right)(\veg r', t)dS'\,,
  \\
      \op T_\text{h}\left(\veg f\right)\left(\veg r,t\right)=\n\left(\veg r\right) \times \nabla \iint_{\veg r' \in \Gamma}\left(\op G_{\veg r}*_t\nabla'\cdot \veg f\right)\left(\veg r',t\right)dS'\,,
\end{gather}
where $c_0$ is the speed of light in the background. The temporal convolution product $*_t$ and the temporal Green function $\op G$ are defined as
\begin{gather}
\left(f*_t g\right)\left(t\right)=\int_{-\infty}^{\infty}f\left(\tau \right)g\left(t-\tau \right)d\tau\, ,\\
\op G_{\veg r}\left(\veg r',t\right)=\frac{\delta\left(t-\frac{|\veg r-\veg r'|}{c_0}\right)}{4\pi |\veg r-\veg r'|}\,,
\end{gather}
with $\delta$ the time Dirac delta.
\subsection{Marching-On-In-Time with Convolution Quadratures}
\label{subsection:CQ_explication}
Let $\Theta$ be a placeholder for any of the integral operators previously presented and let $\veg k(\veg r,t)$ be a causal function (${\forall t<0}$, $\veg k(\veg r,t)=\veg 0$). With these notations, most time domain integral equation take the form
\begin{equation}
    \Theta \left(\veg f_\mr{c}\right)\left(\veg r,t\right)=\veg k\left(\veg r,t\right)\,,
    \label{eq:continuous}
\end{equation}
where $\veg f_\mr{c}$ is the solution to be solved for. The first step of the Marching-On-in-Time solution scheme with convolution quadratures is to apply the boundary element method \cite{harrington1996field, gibson2021method, rao1982electromagnetic} as spatial discretization. Assuming separability between the space and time variables, the unknown function $\veg f_\mr{c}$ is expanded as a linear combination of $N_s$ spatial basis functions such that \cite{banjai2022integral}
\begin{equation}
    \veg f_\mr{c}\left(\veg r,t\right)\approx\sum_{n=1}^{N_s}[\vec f_\Gamma]_n\left(t\right) \veg f^\text{src}_n\left(\veg r\right)\,,
    \label{eq=approx1}
\end{equation}
where $\{\veg f^{\mr{src}}\}$ are the source spatial basis functions and their associated time coefficients are stored in the vector $\vec f_\Gamma(t)$. Then, the equation \eqref{eq:continuous} is tested by the spatial basis functions $\{\veg f^\text{tst}\}$ leading to the time-dependent matrix system
\begin{equation}
    \left(\mat \theta *\vec f\right)\left(t\right)=\vec k\left(t\right)
    \label{eq:time_semi}\,,
\end{equation}
where for $n$ and $m$ in $\llbracket 1,N_s \rrbracket$, we have,
\begin{equation}
    \begin{split}
        \left(\mat \theta *\vec f\right)_m\left(t\right) & =\sum_{n=1}^{N_s}\left( \mat \theta_{m,n}*_t f_n\right)\left(t\right)\,,
        \\
        \begin{bmatrix}
            \mat \theta
        \end{bmatrix}_{m,n}\left(t\right) & = \langle \veg f^\text{tst}_{m},\Theta\left(\delta \veg f^\text{src}_{n}\right)\rangle_\Gamma \left(t\right)\,,
        \\
        \begin{bmatrix}
            \vec k_\Gamma        \end{bmatrix}_m\left(t\right)&=\langle \veg f^\text{tst}_{m},\veg k\rangle_\Gamma \left(t\right)\,,
    \end{split}
\end{equation}
with 
\begin{equation}
    \langle \veg f,\veg g\rangle_\Gamma=\iint_{\veg r \in \Gamma}\veg f\left(\veg r\right) \cdot \veg g\left(\veg r\right) dS\,.
\end{equation}
The second step is the discretization in time with the convolution quadrature method \cite{lubich1988convolutionI,lubich1988convolutionII,wang2011implicit,banjai2022integral,lubich1993runge,banjai2012runge}. 
First, the system \eqref{eq:time_semi} must be transformed in the Laplace domain \cite{schwartz2008mathematics} and we denote $\mat \theta_\mathcal{L}$, $\vec f_\mathcal{L}$ and $\vec k_\mathcal{L}$ the Laplace transform of $\mat \theta$, $\vec f_\Gamma$ and $\vec k_\Gamma$. The system \eqref{eq:time_semi} is then equivalent to
\begin{equation}
   \mathcal{L}\{t \mapsto \left(\mat \theta*\vec f\right) \, \left(t\right) \}(s)=\mat \theta_\mathcal{L}\left(s\right) \vec f_\mathcal{L}\left(s\right)=\vec k_\mathcal{L}\left(s\right)\,.
    \label{eq::Laplace}
\end{equation}
Then, a representation on the Z-domain discretizes the system \eqref{eq::Laplace}. The Laplace parameter, in the operator $\mat \theta_\mathcal{L}$, is replaced by the matrix-valued parameter
\begin{equation}
     \mat s_\mr{cq}\left(z\right)=\frac{1}{\Delta t}\left(\mat A+\frac{\vec 1_p\vec b^T}{z-1}\right)^{-1}\,,
     \label{eq::sCQ}
\end{equation}
diagonalizable for the considered $z$ values, with the following eigenvalue decomposition $\mat s_\mr{cq}(z)=\mat Q(z) \mathbf{\Lambda}(z)\mat Q^{-1}(z)$. The elements of the diagonal matrix $\mathbf{\Lambda}(z)$ are the eigenvalues of $\mat s_\mr{cq}(z)$ and the columns of $\mat Q(z)$ are their associated eigenvectors. The time step size is denoted $\Delta t$ and the matrix $\mat A$ and the vectors $\vec 1_p$, $\vec c$, $\vec b$ of size $p$ are determined by the implicit scheme used \cite{dely2022convolution}. The discretized Z-domain operator $\mat \theta_\mathcal{Z}$ is defined such that for any for $k$ and $l$ in $\llbracket 1,p \rrbracket$ and for $m$ and $n$ in $\llbracket 1,N_s \rrbracket$
\begin{equation}
    \begin{bmatrix}
            \mat \theta_\mathcal{Z}
     \end{bmatrix}_{\Phi_{m,k}, \Phi_{n,l}}\left(z\right)
     =\begin{bmatrix}\mat Q(z)  \mat \theta_{m,n}^{\mathbf{\Lambda}}(z) \mat Q^{-1}(z)\end{bmatrix}_{k,l}\, ,
     \label{eq::Z_operator}
\end{equation}
where $ \Phi_{\alpha,\beta}=p(\alpha-1)+\beta$ is an appropriate indexing function and $\mat \theta_{m,n}^{\mathbf{\Lambda}}(z)$ is a diagonal matrix defined as
\begin{equation}
    \begin{bmatrix}
        \mat \theta_{m,n}^{\mathbf{\Lambda}}
    \end{bmatrix}_{k,k}(z)=
             \begin{bmatrix}\mat \theta_\mathcal{L}({\mathbf{\Lambda}_{k,k}})\end{bmatrix}_{m,n}(z)\,.
\end{equation}
The vectors $\vec f_\mathcal{Z}$ and $\vec k_\mathcal{Z}$ are the Z-domain representation \cite{dely2022convolution} of the respective time-discretized vectors
\begin{equation}
    \begin{split}
       & \begin{bmatrix}
            \vec k_i
        \end{bmatrix}_{\Phi_{m,k}}=\begin{bmatrix} \vec k_\Gamma \end{bmatrix}_m\left( \Delta t \left(i  +[\vec c]_k \right)\right) \,,\\
       & \begin{bmatrix}
            \vec f_i
        \end{bmatrix}_{\Phi_{n,k}}=\begin{bmatrix} \vec f_\Gamma \end{bmatrix}_n \left( \Delta t \left(i  +[\vec c]_k \right)\right) \,,
    \end{split}
\end{equation}
yielding to the following discretization of the system \eqref{eq::Laplace}
\begin{equation}
   \mat \theta_\mathcal{Z}\left(z\right) \vec f_\mathcal{Z} \left(z\right) =\vec k_\mathcal{Z}\left(z\right)\, .
   \label{eq::z_system}
\end{equation}
Finally, the equivalent time discretized system of \eqref{eq:time_semi} is obtained by applying the inverse Z-transform on \eqref{eq::z_system}
\begin{equation}
\mathcal{Z}^{-1}\{\mat \theta_\mathcal{Z}(z)\vec f_\mathcal{Z}(z)\}_i=[\mat Z_{\mat \theta}*_s \veg f]_i= \sum_{j=0}^{i}\mat Z_{\mat \theta,j}\veg f_{i-j} =\vec k_i\, ,
\label{eq:CQMOM}
\end{equation}
where $\mat Z_{\mat \theta,i}=\mathcal{Z}^{-1}\{z \rightarrow \mat \theta_\mathcal{Z}(z)\}_i$ are the time domain interaction matrices and $*_s$ is the sequence convolution product. The system sequence \eqref{eq:CQMOM} is rewritten in the following Marching-On-In-Time that can be solved for $\vec f_i$ 
\begin{equation}
    \quad \mat Z_{\mat \theta,0}\vec f_i =\vec k_i-\sum_{j=1}^{i}\mat Z_{\mat \theta,j}\vec f_{i-j}\, .
    \label{eq:MOT_base}
\end{equation}
\subsection{Classic Integral Marching-On-In-Times}
In this subsection, the discretization scheme described above will be applied to the specific case of the EFIE. The Rao-Wilton-Glisson (RWG) basis functions $\{ \veg f^{\mr{rwg}}_n \}_{N_s}$ \cite{rao1982electromagnetic,raviart1977mixed} are used to expand the current density as
\begin{equation}
    \veg j_\Gamma \left(\veg r,t\right) \approx \sum_{n=1}^{N_s} [ \vec j_\Gamma]_n \left(t\right)\veg f^{\mr{rwg}}_n\left(\veg r\right)\, ,
\end{equation}
where the current coefficients are gathered in an unknown vector function of time $ \vec j_\Gamma (t)$.
The EFIE is then tested with rotated RWG basis functions $\{\n \times\veg f^{\mr{rwg}}_n\}_{N_s}$, leading to the following Marching-On-In-Time
\begin{equation}
\mat Z_{\mat T, 0}\vec j_i =-\eta_0^{-1} \vec e^\mr{inc}_i - \sum_{j=1}^{i}\mat Z_{\mat T, j}\vec j_{i- j}\,,
 \label{eq:MOTEFIEunbounded}
\end{equation}
where the vector sequences $\vec j_i$ and $\vec e^\mr{inc}_i$, and the time domain interaction matrices $\mat Z_{\mat T,i}$ are respectively generated by the convolution quadrature method described in Subsection~\ref{subsection:CQ_explication} of $\vec j_\Gamma (t) $ and the following space-discretized vector and matrix
\begin{gather}
    \begin{bmatrix}
    \vec e^\mr{inc}_\Gamma \end{bmatrix}_m(t)
= \langle \n \times \veg f^\text{rwg}_{m}, \n  \times \veg e^\mr{inc} \rangle_\Gamma \left( t \right)\, ,\\
    \begin{bmatrix}\mat T\end{bmatrix}_{m,n}(t)=\langle \n \times \veg f^\mr{rwg}_{m},\op T \left(\delta \veg f^\mr{rwg}_{n}\right) \rangle_\Gamma \left( t \right)\, .
\end{gather}
However, the time integral contribution of this operator $\op T$ involves an unbounded number of non-vanishing matrices $\mat Z_{\mat T,i}$ \eqref{eq:MOTEFIEunbounded}, leading to a prohibitive quadratic complexity with the number of time steps \cite{beghein2013temporal}. Historically, the time differentiated formulation is preferred because it is not afflicted by this drawback
\cite{wang2011implicit,dely2019large}, and leads to the following MOT \cite{dely2019large}
\begin{equation}
\mat Z_{\dot{\mat T},0}\vec j_i =-\eta_0^{-1} \dot{\vec e}^\mr{inc}_i - \sum_{j=1}^{N_\mr{conv}}\mat Z_{\dot{\mat T},j}\vec j_{i- j}\,,
 \label{eq:MOTTDEFIE}
\end{equation}
where $\dot{\vec e}^\mr{inc}_n$ and $\dot{\mat Z}_n$ are respectively the time domain vectors and interaction matrices generated by the convolution quadrature method described in Subsection~\ref{subsection:CQ_explication} of 
\begin{gather}
 \begin{bmatrix}
     \dot{\vec e}^\mr{inc}_\Gamma
 \end{bmatrix}_m(t)=\langle \n \times \veg f^\text{rwg}_{m}, \n  \times \frac{\partial}{\partial t}\veg e^\mr{inc} \rangle_\Gamma\,,\\
    \begin{bmatrix}
      \dot{\mat T}   
     \end{bmatrix}_{m,n}(t)=\langle \n \times \veg f^\mr{rwg}_{m},\frac{\partial}{\partial t}\op T \left(\delta \veg f^\mr{rwg}_{n}\right) \rangle_\Gamma \, .
    \end{gather}
\subsection{EFIE DC instability}
\begin{figure}[ht]
    \includegraphics[scale=0.9]{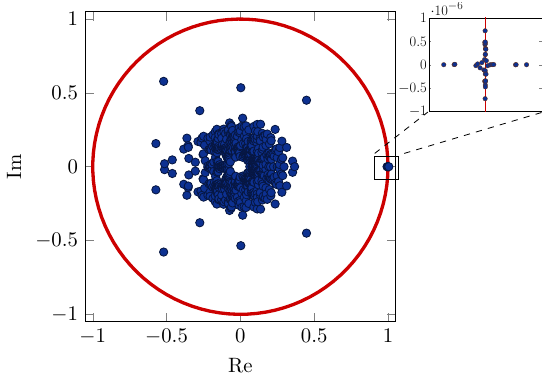}
    \caption{Polynomial eigenvalues of the TD-EFIE MOT scheme on the unit sphere, $N_s=270$ and $\Delta t=\SI{3}{\nano\second}$}
    \label{fig:sphere_eigenvalues}
\end{figure}

The electric field integral operator suffers from the DC instability: since for all constant-in-time solenoidal current $\veg j_\mr{cs}$ we have $\nabla \cdot \vec j_\mr{cs}=0$ and $\frac{\partial}{\partial t}\veg j_\mr{cs}=\mathbf{0}$, we can conclude that
\begin{equation}
    \op T \left(\veg j_\mr{cs}\right)=\mathbf{0}\, .
\end{equation}
Therefore, the EFIE solution is only determined up to a constant solenoidal current \cite{dely2019large}. Its time differentiated counterpart inherits these drawbacks and amplifies the DC instability by further adding linear in time solenoidal currents to the nullspace. This latter deteriorates the late time simulation in which spurious currents grow exponentially in the operator nullspace \cite{van2022role}. This behaviour is predicted by the polynomial eigenvalues analysis of the MOT: a stable MOT has all its eigenvalues inside the unit circle in the complex plane while a MOT that suffers from the DC instability has some eigenvalues that cluster around $1$ \cite{walker2002stability}. The eigenvalue distribution of the time differentiated EFIE MOT is represented in \figref{fig:sphere_eigenvalues} in which such a cluster is clearly visible around $1$.
\subsection{Quasi-Helmholtz projectors}
Previous works show that the electric field integral equation discretized in space using RWG basis functions can be stabilized by the quasi-Helmholtz projectors \cite{ dely2019large, beghein2015dc}. These projectors are formed from the star-to-rwg transformation matrix, denoted $\mat \Sigma$ and defined in \cite{adrian2021electromagnetic}, which maps the discretized current into the non-solenoidal contributions \cite{9580445,nedelec2001acoustic}. The quasi-Helmholtz projectors on the non-solenoidal space and its complementary (the one on solenoidal/quasi-harmonic space) are respectively
\begin{equation}
    \mat P^{\Sigma}=\left(\mat \Sigma \left(\mat \Sigma^T\mat \Sigma \right)^{+}\mat \Sigma^T\right)\quad \text{and}\quad \mat P^{\Lambda H}=\mat I-\mat P^{\Sigma}\, ,
\end{equation}
where $\text{ }^+$ denotes the Moore–Penrose pseudoinverse \cite{adrian2021electromagnetic}. 
\section{Calderón preconditioning of the EFIE}
\label{section:CP_EFIE}
EFIE formulations based on the quasi-Helmholtz projectors cure the DC instability and the conditioning at large time steps. However, these formulations still suffer from a dense discretization breakdown. One appealing strategy could be to apply standard preconditioning schemes to the Marching-On-In-Time matrices directly to cure the matrix conditioning issues. However, the solution currents would remain unaltered and subject to DC instabilities as the original scheme. This is why the preconditioning has to be performed on the continuous equations to build a new operator without nullspace and then discretize the formulation to obtain a well-conditioned scheme. In this part, Calderón preconditioning strategies are proposed to cure the DC instability and the conditioning breakdowns.
\subsection{A Convolution Quadrature Calderón time-domain EFIE}
\begin{figure}[ht]
\centering
    \includegraphics[scale=0.9]{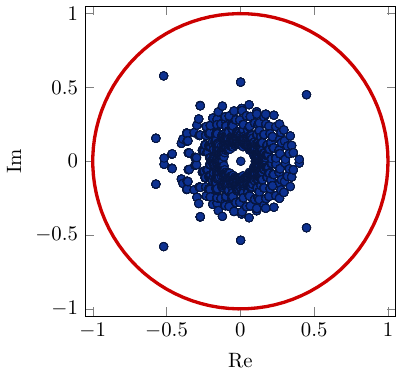}
    \caption{Polynomial eigenvalues of the Calderón EFIE MOT scheme on the unit sphere, $N_s=270$ and $\Delta t=\SI{3}{\nano\second}$.}
    \label{fig:cal_eigen_sphere}
\end{figure}

\begin{figure}[ht]
\centering
    \includegraphics[scale=0.9]{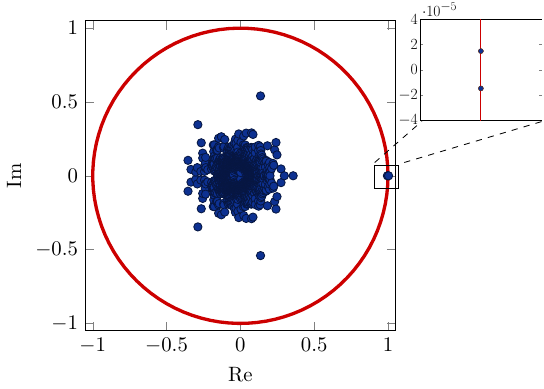}
    \caption{Polynomial eigenvalues of the Calderón EFIE MOT scheme on a torus with a inner and outer radii respectively equal to $\SI{0.2}{\metre}$ and $\SI{0.5}{\metre}$, $N_s=387$ and $\Delta t=\SI{3}{\nano\second}$. Near $1$, four eigenvalues are clustered, superposed two by two on this figure.}
    \label{fig:cal_eigen_torus}
\end{figure}

Calderón preconditioners are based on the Calderón identity \cite{hsiao1997mathematical,cools2009time} 
\begin{equation}
    \op T^2=\op T  \circ \op T= -\op I /4+ \op K^2\, ,
    \label{eq:calIdendity}
\end{equation}
where $\circ$ is the composition operator, $\op I$ is the identity operator and $\op K$ is defined as
\begin{equation}
    \op K\left(\veg f\right)\left(\veg r,t\right)=\n\left(\veg r\right) \times \nabla \times \iint_{\veg r' \in \Gamma}\left(\op G_{\veg r}*_t \veg f\right) (\veg r',t)dS'\,.
\end{equation}
The operator $-\op I /4+ \op K^2$ is a well-behaved operator for increasing discretization densities. As a consequence, with a proper discretization, $\op T^2$ is well-conditioned for large time steps and dense meshes for simply connected structures \cite{adrian2021electromagnetic}. In practice, a discretization of $\op T^2$ is used in which the right EFIE operator is discretized with RWG basis functions and the left preconditioner is discretized with Buffa-Christiansen (BC) basis functions $\left(\veg f^{\mr{bc}}_n\right)_{N_s}$ \cite{andriulli2012analysis, buffa2007dual,cools2011accurate}
\begin{equation}
\begin{bmatrix}\mathbb{T}\end{bmatrix}_{m,n} \left(t\right)=\langle\n \times \veg f^\mr{bc}_{m},\op T \left(\delta \veg f^\mr{bc}_{n} \right)\rangle_\Gamma (t)\,.
        \label{eq:BC_EFIO}
\end{equation}
The preconditioning leads to the following space-discretized formulation
\begin{equation}
   \left(\bm{\mathbb{T}} \mat G_m^{-1} * \left(\mat T * \vec j_\Gamma\right)\right) \left(t \right)=-\eta_0^{-1}\left(\bm{\mathbb{T}} \mat G_m^{-1} *\vec e^\mr{inc}_\Gamma\right)\left(t \right)\, ,
\end{equation}
where the matrix $\mat G_m$ is the mixed gram matrix linking the the two discretizations
\begin{equation}
    \begin{bmatrix}
        \mat G_m
    \end{bmatrix}_{m,n}=\langle \n \times \veg f^{\mr{rwg}}_m,\veg f^{\mr{bc}}_n\rangle_\Gamma\,.
\end{equation}
Then, the convolution quadrature leads to the MOT scheme
\begin{equation}
\begin{split}
    \begin{bmatrix}
    \mat Z_\mathbb{T}\widetilde{\mat G}_m^{-1}*_s \mat Z_{\mat T}
\end{bmatrix}_{0} \vec j_{i} &=-\eta_0^{-1} \begin{bmatrix}\mat Z_\mathbb{T}\widetilde{\mat G}_m^{-1} *_s\vec e^\mr{inc}\end{bmatrix}_i
 \\
 &-\sum_{j=1}^{i}\begin{bmatrix}
    \mat Z_\mathbb{T}\widetilde{\mat G}_m^{-1}*_s \mat Z_{\mat T}
\end{bmatrix}_{j}\vec j_{i- j}\,,
\end{split}
  \label{CMP}
\end{equation}
where $\mat Z_{\mathbb{T},i}$ are the time domain interaction matrices of the space-discretized operator $\mathbb{T}\left(t\right)$ \eqref{eq:BC_EFIO} generated by the convolution quadrature method described in Subsection~\ref{subsection:CQ_explication}, the sequence convolution quadrature product $*_s$ is the discretization of the space-discretized temporal convolution product $*$ and 
\begin{equation}
    \widetilde{\mat G}_m=\mat G_m \otimes\mat I_p \, .  
\end{equation}
The Kronecker product $\otimes\mat I_p$ is required to match with the convolution quadrature method where $\mat I_p$ is the identity matrix of size $p$.
Unfortunately, the MOT in \eqref{CMP}, involves operators with temporal integrations leading to a time consuming MOT. A more favorable scheme can be obtained by noticing the following commutative properties
\begin{gather}
    \frac{\partial}{\partial t}\op T_{\text{s}/\text{h}}\left(\veg f\right)\left(\veg r,t\right)=\op T_{\text{s}/\text{h}}\left(\frac{\partial}{\partial t}\veg f\right)\left(\veg r,t\right), \\
    \int_{-\infty}^{t}\op T_{\text{s}/\text{h}}\left(\veg f\right)\left(\veg r,t'\right) dt'=\op T_{\text{s}/\text{h}}\left(\int_{-\infty}^{.}\veg f\right)\left(\veg r,t\right)\, ,\label{eq::int_com}
\end{gather}
and the cancellation property $\op T_\text{h}^2=\mathbf{0}$ \cite{andriulli2008multiplicative}, we have
\begin{align}
\begin{split}
    \op T^2&=\left(-\frac{1}{c_0}\frac{\partial}{\partial t}\op T_\text{s} + c_0\int_{-\infty}^{\cdot} \hspace{-7px} \op T_\text{h}\right)\circ\left(-\frac{1}{c_0}\frac{\partial}{\partial t}\op T_\text{s} + c_0\int_{-\infty}^{\cdot} \hspace{-7px} \op T_\text{h}\right)\\
    &=c_0^{-2}\frac{\partial^2}{\partial t^2}\op T_\text{s}^2-  \op T_\text{s}\circ \op T_\text{h}-\op T_\text{h}\circ \op T_\text{s} \,.
\end{split}
\label{eq:Calderon_op}
\end{align}
This is advantageous since, besides not involving any time integration contribution, the operator $c_0^{-2}\frac{\partial^2}{\partial t^2}\op T_\text{s}^2- \op T_\text{s}\op T_\text{h}-\op T_\text{h}\op T_\text{s}$ has no nullspace for simply connected geometries leading to a DC-stable discretization (``dottrick TDEFIE'') \cite{andriulli2008dottrick}. By extending the
previous notations on $\op T_\text{s}$ and $\op T_\text{h}$
\begin{equation}
    \begin{split}
    &\begin{bmatrix}\mat T_{\text{s}/\text{h}}\end{bmatrix}_{m,n}(t)=\langle \n \times \veg f^\mr{rwg}_{m},\op T_{\text{s}/\text{h}} \left(\delta \veg f^\mr{rwg}_{n}\right) \rangle_\Gamma \left( t \right)\, ,
    \\
&\begin{bmatrix}\mathbb{T}_{\text{s}/\text{h}}\end{bmatrix}_{m,n}\left(t\right)  =\langle\n \times \veg f^\mr{bc}_{m},\op T_{\text{s}/\text{h}} \left(\delta \veg f^\mr{bc}_{n} \right)\rangle_\Gamma \left( t \right)\, ,
    \end{split}
\end{equation}
the proposed space-discretized operator is denoted 
\begin{equation}
     \mr{T}_{c}=c_0^{-2}\frac{\partial^2}{\partial t^2}\bm{\mathbb{T}}_{\text{s}} \mat G_m^{-1} *\mat T_{\text{s}} -\bm{\mathbb{T}}_{\text{s}} \mat G_m^{-1}* \mat T_{\text{h}}
     -\bm{\mathbb{T}}_{\text{h}} \mat G_m^{-1}* \mat T_{\text{s}}\,,
     \label{eq:Calderon_operator}
\end{equation}
yielding to the following space-discretized formulation
\begin{equation}
\left( \mr{T}_{c}* \vec j_\Gamma \right)(t)=-\eta_0^{-1}\left(\bm{\mathbb{T}}\mat G_m^{-1}* \vec e^\mr{inc}_\Gamma\right)(t)\, .
\label{eq::Laplace_form}
\end{equation}
The right-hand side operator still involves a temporal integration in \eqref{eq::Laplace_form}. However, given the commutative properties \eqref{eq::int_com}, the temporal integral on the scalar potential $\op T_h$ is evaluated with the incident field
\begin{equation}
 c_0 \left(\int_{-\infty}^{\cdot} \hspace{-7px}\bm{\mathbb{T}}_{\text{h}} \mat G_m^{-1}* \vec e^\mr{inc}_\Gamma\right)(t)= c_0 \left(\bm{\mathbb{T}}_{\text{h}} \mat G_m^{-1}* \vec e^\mr{prim}_\Gamma\right)(t)\,,
\end{equation}
where
\begin{equation}
\begin{split}
    &\begin{bmatrix}
    \vec e^\mr{prim}_\Gamma \end{bmatrix}_m(t)
= \langle \n \times \veg f^\text{rwg}_{m}, \int_{-\infty}^{t} \n \times \veg e^\mr{inc}(t') dt'\rangle_\Gamma \,,\\
    &\begin{bmatrix} \vec e^\mr{prim}_i\end{bmatrix}_{\Phi_{m,k}}=
    \begin{bmatrix} \vec e^\mr{prim}_\Gamma \end{bmatrix}_m \left(\Delta t \left(i+[\vec c]_k \right)\right)\,.
\end{split}
\end{equation}
Therefore, the previous MOT is rewritten as 
\begin{equation}
\mat Z_{\mr{T}_{\text{c}},0}
\vec j_{i}=- \eta_0^{-1}\sum_{j=0}^{N_\mr{conv}}\left( \mat Z_{\dot{\mr{T}}_\text{s},j}\vec e^\mr{inc}_{i- j}+\mat Z_{\mr{T}_\text{h},j}\vec e^\mr{prim}_{i- j}\right)
 -\sum_{j=1}^{N_\mr{conv}}
Z_{\mr{T}_{\text{c}},j}\vec j_{i- j}\,,
 \label{eq:MOTCPEFIE}
\end{equation}
where the time domain interaction matrices $\mat Z_{\mr{T}_{\text{c}},i}$, $ \mat Z_{\dot{\mr{T}}_\text{s},i}$ and $\mat Z_{\mr{T}_\text{h},i}$ are respectively generated by the convolution quadrature method described in Subsection~\ref{subsection:CQ_explication} of the space-discretized operators $\mr{T}_\text{c}$, $\dot{\mr{T}}_\text{s}=-c_0^{-1}\frac{\partial}{\partial t}\bm{\mathbb{T}}_{\text{s}} \mat G_m^{-1}$ and $\mr{T}_\text{h}=c_0\bm{\mathbb{T}}_{\text{h}} \mat G_m^{-1}$. 
As in \eqref{CMP}, the interaction matrix sequence $\mat Z_{\mr{T}_\text{c},i}$ involves computationally expensive sequence convolution products $*_s$, however, the convolution quadrature method allows the substitution of the sequence convolution products $*_s$ by matrix multiplications in the Z-domain, that can be evaluated at a lesser cost. By extending the notations of the convolution quadrature described in Subsection~\ref{subsection:CQ_explication} on the space-discretized operators $\mat T_\text{s/h}$ and $\mathbb{T}_\text{s/h}$ and by using the Z-domain properties, the matrix sequence $\mat Z_{\mr{T}_\text{c},i}$ is equal to
\begin{equation}
\begin{split}
    \mat Z_{\mr{T}_\text{c},i}& = c_0^{-2}\mathcal{Z}^{-1}\{\overline{\mat s}_\mr{cq}^2\mathbb{T}_{\text{s},\mathcal{Z}} \widetilde{\mat G}_m \mat T_{\text{s},\mathcal{Z}} \} \\&  -\mathcal{Z}^{-1}\{\mathbb{T}_{\text{h},\mathcal{Z}}  \widetilde{\mat G}_m \mat T_{\text{s},\mathcal{Z}} +\mathbb{T}_{\text{s},\mathcal{Z}}  \widetilde{\mat G}_m \mat T_{\text{h},\mathcal{Z}} \} \, ,
\end{split}
\end{equation}
where the matrix $\overline{\mat s}_\mr{cq}(z)=\mat I_{N_s}\otimes\mat s_\mr{cq}(z)$ is the Z-discretization of the time derivative and $\mat I_{N_s}$ is the identity matrix of size $N_s$.
The formulation \eqref{eq:MOTCPEFIE} is a good candidate to obtain a stable current solution, however, the proposed operator $\op T^2$ has static nullspaces for multiply connected geometries \cite{cools2009nullspaces}. As such, \eqref{eq:Calderon_operator} is still subject to DC-instabilities for multiply connected geometries. The polynomial eigenvalue analysis on a sphere and on a torus (respectively \figref{fig:cal_eigen_sphere} and \figref{fig:cal_eigen_torus}) illustrate this phenomenon. While all the eigenvalues cluster in $0$ in the spherical case, an analysis on a torus highlights four eigenvalues of this MOT clustered around $1$, corresponding to the four constant regime solutions \cite{cools2009nullspaces}.
\subsection{A Convolution Quadrature Calderón time-domain EFIE regularized with quasi-Helmholtz projectors}
The previous Calderón formulation is perfectly adapted to simply connected geometries, ensuring that the new operator has no nullspace. However, on multiply-connected geometries, the harmonic subspace is non-empty, thus enlarging the nullspace of $\op T$ which is a new source of DC instability in \eqref{eq:MOTCPEFIE}\cite{cools2009nullspaces}. The discretized EFIE operators can be regularized using the quasi-Helmholtz projectors to address this issue. Because the regularization is based on projectors, it does not compromise the $h$-refinement regularizing effect of the original Calderón scheme. The regularized EFIE space-discretized operators are 
\begin{gather}
    \mat T^\mr{reg}=\left(\frac{c_0}{a}\int_{-\infty}^{\cdot} \hspace{-7px}\mat P^{\Lambda H}+\mat P^{\Sigma}\right)\mat *\mat T*\left(\mat P^{\Lambda H}+\frac{a}{c_0} \frac{\partial}{\partial t}\mat P^{\Sigma} \right)\, ,
    \label{eq:qh-Laplace}
    \\
    \mathbb{T}^\mr{reg} =\left(\frac{c_0}{a}\int_{-\infty}^{\cdot} \hspace{-7px}\mathbb{P}^{\Sigma H}+\mathbb{P}^{\Lambda}\right)*\mathbb{T}*\left(\mathbb{P}^{\Sigma H}+\frac{a}{c_0}\frac{\partial}{\partial t}\mathbb{P}^{\Lambda}\right)\, ,
    \label{eq::regBC}
\end{gather}
where the BC quasi-Helmholtz projectors are defined with the loop-to-RWG transformation matrix $\mat \Lambda$ \cite{adrian2021electromagnetic} such that
\begin{equation}
    \mathbb{P}^{\Lambda}=\left(\mat \Lambda\left(\mat \Lambda^T\mat \Lambda \right)^{+}\mat \Lambda^T\right) \quad \text{and}\quad \mathbb{P}^{\Sigma H}=\mat I-\mathbb{P}^{\Lambda}\,,
\end{equation}
and where the scaling $\frac{c_0}{a}$ with $a$ defined as the maximal diameter of the scatterer, ensures consistent dimensionality and helps reduce the conditioning further. This application of the projectors is equivalent to differentiating the non-solenoidal contributions on the left and in time integrating the solenoidal contributions on the right of each EFIE operator \cite{dely2019large}.
The regularized Calderón operator in the space-discretized time domain is
\begin{equation}
\mr{T}^\mr{reg}_{\text{c}}= \mathbb{T}^{\mr{reg}} \mat G_m^{-1} * \mat T^{\mr{reg}}\, .
\end{equation}
At first sight, \eqref{eq:qh-Laplace} and \eqref{eq::regBC} seem to involve unpractical temporal integrals. However, the problematic contributions in the regularized EFIE operator {$\mathbf{T}^{\mr{reg}}$ will vanish, since $\mat P^{\Lambda H}\mat T_{\text{h}}=\mat T_{\text{h}} \mat P^{\Lambda H} = \mathbf{0}$ and $\mat P^{\Sigma}\mat T_{\text{h}}\mat P^{\Sigma} =\mat T_{\text{h}}$, and we have 
\begin{equation}
\begin{split}
    \mat T^{\mr{reg}}&=-a^{-1}\mat P^{\Lambda H}\mat T_{\text{s}} \mat P^{\Lambda H} - c_0^{-1} \mat P^{\Lambda H}\frac{\partial}{\partial t}\mat T_{\text{s}} \mat P^\Sigma
     \\& -
     c_0^{-1} \mat P^\Sigma \frac{\partial}{\partial t}\mat T_{\text{s}} \mat P^{\Lambda H} -\frac{a}{c_0^2} \mat P^\Sigma \frac{\partial^2}{\partial t^2}\mat T_{\text{s}} \mat P^\Sigma+ a\mat T_{\text{h}}\, .
\end{split}
\end{equation}
Similarly, the dual EFIE operator simplifies as
\begin{equation}
\begin{split}
       \mathbb{T}^{\mr{reg}}&= -a^{-1}\mathbb{P}^{\Sigma H}\mathbb{T}_{\text{s}} \mathbb{P}^{\Sigma H} - c_0^{-1} \mathbb{P}^{\Sigma H} \frac{\partial}{\partial t}\mathbb{T}_{\text{s}} \mathbb{P}^{\Lambda}
       \\& -  c_0^{-1} \mathbb{P}^{\Lambda} \frac{\partial}{\partial t}\mathbb{T}_{\text{s}} \mathbb{P}^{\Sigma H} -\frac{a}{c_0^2}  \mathbb{P}^{\Lambda} \frac{\partial^2}{\partial t^2} \mathbb{T}_{\text{s}} \mathbb{P}^{\Lambda} + a\mathbb{T}_{\text{h}}\, .
\end{split}
\end{equation}
The space-discretized formulation of the regularized Calderón EFIE is
\begin{equation}
    \mr{T}^{\mr{reg}}_\text{c}* \vec y_\Gamma =-\eta_0^{-1}\mathbb{T}^{\mr{reg}} \mat G_m^{-1}*\mat R* \vec e^\mr{inc}_\Gamma \, ,
    \label{eq:RegCPEFIE}
\end{equation}
where
\begin{equation}
    \mat R =\frac{c_0}{a} \int_{-\infty}^{\cdot} \hspace{-7px}\mat P^{\Lambda H}+\mat P^{\Sigma}\quad \text{and}\quad \vec j_\Gamma= \left(\mat P^{\Lambda H}+\frac{a}{c_0} \frac{\partial}{\partial t}\mat P^{\Sigma} \right)\vec y_\Gamma \, .
    \label{eq:regy2j}
\end{equation}
The right-hand side of \eqref{eq:RegCPEFIE} has a temporal integral which is directly evaluated on the incident field to avoid quadratic complexity with the number of time steps. This leads to the following MOT 
\begin{equation}
\begin{split}
    \mat Z_{\mr{T}_\text{c}^{\mr{reg}},0}\mat y_{i}&=-\eta_0^{-1}\sum_{j=0}^{N_\mr{conv}}\mat Z_{\mr{T}^{\mr{reg}},j} \left(\widetilde{\mat P}^{\Sigma}\vec e^\mr{inc}_{i- j} +\frac{c_0}{a} \widetilde{\mat P}^{\Lambda H} \vec e^\mr{prim}_{i- j} \right)\\
    &-\sum_{j=1}^{N_\mr{conv}} \mat Z_{\mr{T}^\mr{reg}_\text{c},j}\vec y_{i- j}\,,
\end{split}
\label{eq:MOTqHCPEFIE}
\end{equation}
where $\widetilde{\mat P}^{\Lambda H}=\mat P^{\Lambda H} \otimes I_p$, $\widetilde{\mat P}^{\Sigma}=\mat P^{\Sigma} \otimes I_p$, the vector sequence $\vec y_n$ is the time discretization of $\vec y_\Gamma (t)$, and the time domain interaction matrices $\mat Z_{\mr{T}^\mr{reg}_\text{c},i}$ and $ \mat Z_{\mr{T}^\mr{reg},i}$ are respectively generated by the convolution quadrature method described in Subsection~\ref{subsection:CQ_explication} of the space-discretized operators $\mr{T}^{\mr{reg}}_\text{c}$ and $\mr{T}^{\mr{reg}}=\mathbb{T}^{\mr{reg}}\mat G_m^{-1}$. Once the computation of $\vec y_n$ is done, the current $\vec j$ still has to be evaluated. The convolution quadrature discretization of the time derivative is 
\begin{equation}
\begin{split}
     \mat Z_{\frac{\partial}{\partial t},i}&=\mathcal{Z}^{-1}\{z\rightarrow\overline{\mat s}_\mr{cq}(z)\}_i\\&=
    \Delta t^{-1} \overline{\mat A^{-1}}\delta_{i,0}-\Delta t^{-1} \overline{\mat A^{-1}\vec 1_p\vec b^{T}\mat A^{-1}}\delta_{i,1}\, ,
\end{split}
    \label{eq::Zinv_s}
\end{equation}
where $\delta_{i,0}$ is the Kronecker delta, $\overline{\mat A^{-1}}=I_{N_s} \otimes \mat A^{-1}$ and $\overline{\mat A^{-1}\vec 1_p\vec b^{T}\mat A^{-1}} =I_{N_s} \otimes  \mat A^{-1}\vec 1_p\vec b^{T}\mat A^{-1}$ \cite{dely2019large}. Therefore, the current solution is obtained as
\begin{equation}
\begin{split}
    \vec j_i&=\left[\left(\widetilde{\mat P}^{\Lambda H}+\frac{a}{c_0}\mat Z_{\frac{\partial}{\partial t}}\widetilde{\mat P}^{\Sigma}\right)*_s\vec y\right]_i\\
    &=\widetilde{\mat P}^{\Lambda H} \vec y_i + \frac{a}{\Delta t c_0} \widetilde{\mat P}^{\Sigma} \left[ \overline{\mat A^{-1}} \vec y_i - \overline{\mat A^{-1}\vec 1_p\vec b^{T}\mat A^{-1}} \vec y_{i-1}\right]\, .
\end{split}
\end{equation}
\section{Results}
\label{section:results}
To test the effectiveness of the proposed schemes, simulations have been realized with different geometries, excited by a Gaussian pulse plane wave
\begin{gather}
 \veg e^\mr{inc}\left(\veg r,t\right)  =A_0 \exp \Big(-\frac{\big(t-\frac{\hat{\veg k}\cdot \veg r}{c}\big)^2}{2\sigma^2}\Big) \hat{\veg p}\, ,
 \\
 \veg h^\mr{inc}\left(\veg r,t\right)  =  \frac{1}{\eta_0}\hat{\veg k}\times \veg e^\mr{inc}\left(\veg r,t\right)\,, 
\end{gather}
where $\sigma=6/(2\pi f_\mr{bw})$, $\hat{\veg p}=\hat{\veg x}$, $\hat{\veg k}=-\hat{\veg z}$, $A_0=\SI{1}{\volt\per\meter}$ and $f_\mr{bw}$ is the frequency bandwidth. Notice that this frequency bandwidth is proportional to the maximal frequencies excited by the pulse Gaussian. In this work, the Runge-Kutta Radau IIA method of stage 2 is used for all simulations \cite{runge1895numerische,kutta1901beitrag}. The time step size $\Delta t$ has been chosen equal to $(\psi f_\mr{max})^{-1}$ where $f_\mr{max}$ is the upper frequency of the excitation and $\psi=3$ is an oversampling parameter.

\subsection{Canonical geometries}
To illustrate the key properties of the newly proposed schemes, namely the Calderón preconditioned formulation \eqref{eq:MOTCPEFIE} and the and the Calderón preconditioned formulation regularized by the quasi-Helmholtz-projectors \eqref{eq:MOTqHCPEFIE}, they are compared in the case of modelling of canonical scatterers to other formulations present on the literature: the EFIE MOT schemes (MOT EFIE) \eqref{eq:MOTEFIEunbounded}, the time-differentiated one (MOT TD-EFIE) \eqref{eq:MOTTDEFIE}, the formulation regularized by the quasi-Helmholtz-projectors (MOT qH-EFIE) \cite{dely2019large}.
In this subsection, the excitation parameters have been chosen not to excite the first resonant mode of the geometries.

The first set of numerical tests were performed on the unit sphere, discretized with $270$ RWG functions. The intensity of the resulting currents at one point of the geometry are shown in \figref{fig:sphere_current_desity}. As expected, the time differentiated and the non-differentiated EFIEs are the only formulations suffering from DC-instabilities on this simply connected scenario. In addition, the condition number of the matrices to invert for each MOT are presented in \figref{fig:cond_number_dt} and \figref{fig:cond_number_ns} with respect to the time step size $\Delta t$ and the mesh density $h^{-1}$. The standard EFIE formulation and its time-differentiated counterpart suffer from ill-conditioning at large time steps while the stabilized ones remain well-conditioned. Instead, only the Calderón preconditioned formulations presented in this work remain well-conditioned for dense discretizations (\figref{fig:cond_number_ns}).

\begin{figure}[ht]
    \centering
    \includegraphics[scale=0.63]{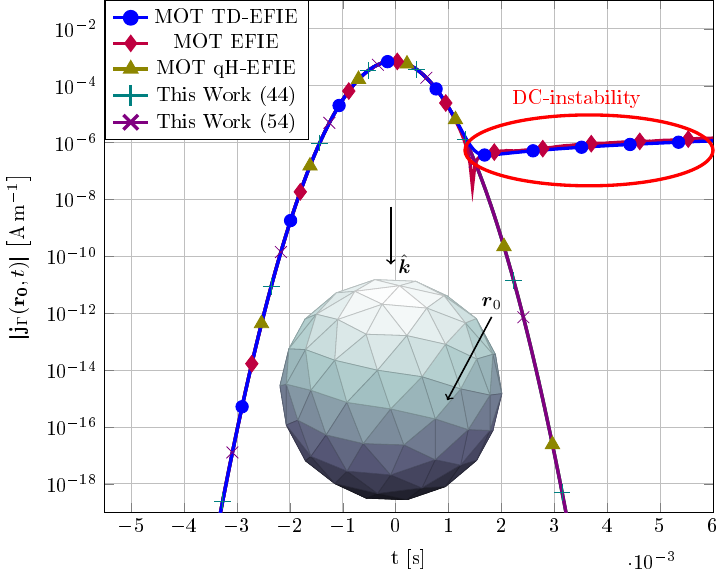}
    \caption{Evolution in time of the current intensity at $\veg r_0=(-0.36,0.89,0.11) \si{\metre}$ of the sphere where $f_{\mr{bw}} =\SI{25}{\kilo\hertz}$ with parameters $N_s=270$ and $\Delta t=\SI{4.5e4}{\nano\second}$.}
    \label{fig:sphere_current_desity}
\end{figure}

\begin{figure}[ht]
    \centering
    \includegraphics[scale=0.63] {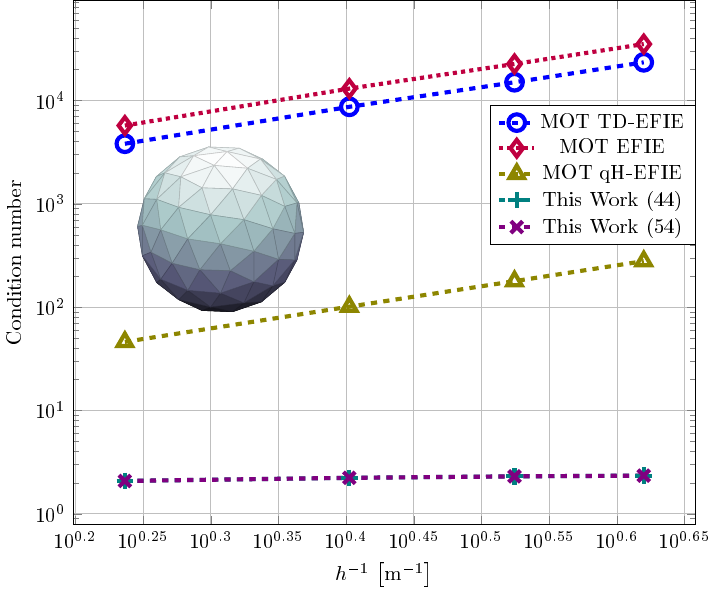}
    \caption{Condition number with respect to the time step size ($N_s=270$) on a sphere of the EFIE MOT schemes.}
    \label{fig:cond_number_dt}
\end{figure}
\begin{figure}[ht]
    \centering
    \includegraphics [scale=0.63]{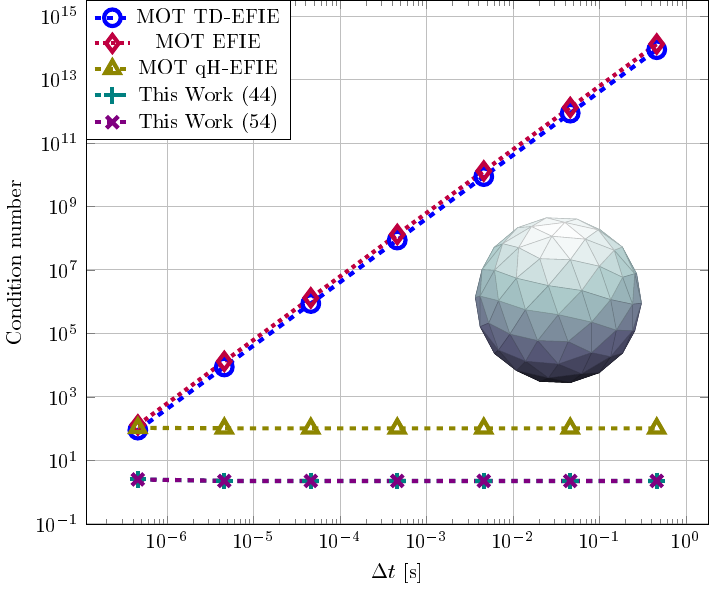}
    \caption {Condition number with respect to the mesh size $h$ ($\Delta t= \SI{45}{\nano\second}$) on a sphere of the EFIE MOT schemes.}
    \label{fig:cond_number_ns}
\end{figure}

The second set of tests focused on the stability of the different formulations when modelling multiply connected scatters, here a torus with inner radius of $\SI{0.2}{\meter}$ and outer radius of $\SI{0.5}{\meter}$. The current densities at the probe point are shown in \figref{fig:torus_current_desity} and the conditioning studies are represented in \figref{fig:torus_cond_number_ns} and \figref{fig:torus_cond_number_dt}. In line with the polynomial eigenvalue analysis of the non-regularized Calderón EFIE formulation (\figref{fig:cal_eigen_torus}), the formulation \eqref{eq:MOTCPEFIE} suffers from DC instability. Moreover, the static nullspace of the continuous operator deteriorates the condition number of the matrix to invert for large time steps (\figref{fig:torus_cond_number_dt}). However, the newly proposed regularized Calderón formulation \eqref{eq:RegCPEFIE} is stable and remains well-conditioned at large time steps and dense meshes for this geometry.
\begin{figure}[ht]
    \includegraphics[scale=0.63]{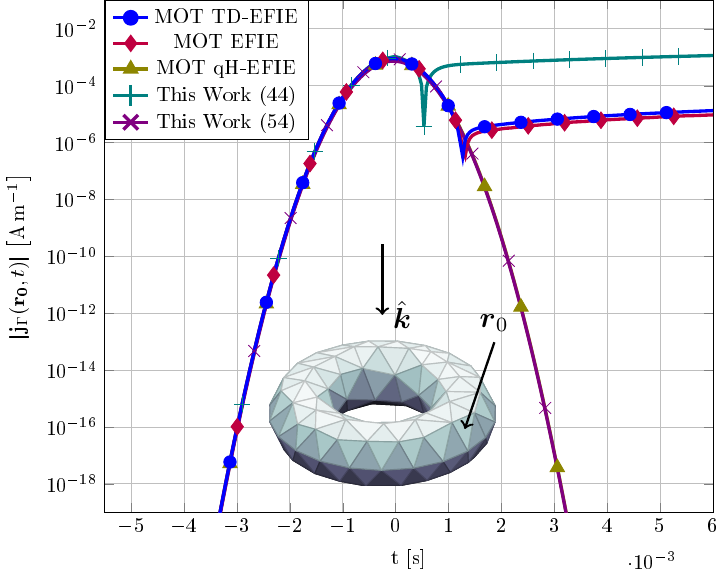}
    \caption{Evolution in time of the current intensity at $\veg r_0=(0.62,-0.13,0.11) \si{\metre}$ of the torus where $f_\mr{bw}=\SI{25}{\kilo\hertz}$ with parameters $N_s=387$ and $\Delta t=\SI{4.5e4}{\nano\second}$.}
    \label{fig:torus_current_desity}
\end{figure}
\begin{figure}[ht]
    \centering
    \includegraphics[scale=0.63]{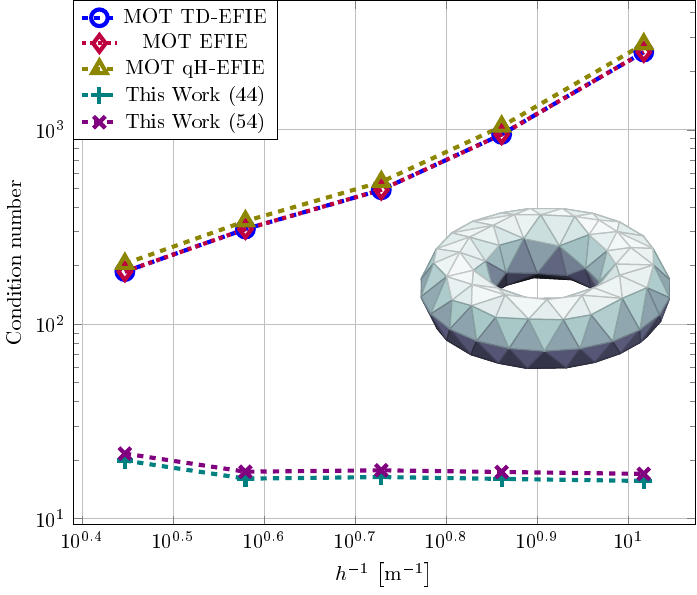}
    \caption{Condition number with respect to the mesh size $h$ ($\Delta t= \SI{4.5}{\nano\second}$) on a torus of the EFIE MOT schemes.}
    \label{fig:torus_cond_number_ns}
\end{figure}
\begin{figure}[ht]
    \centering
    \includegraphics[scale=0.63] {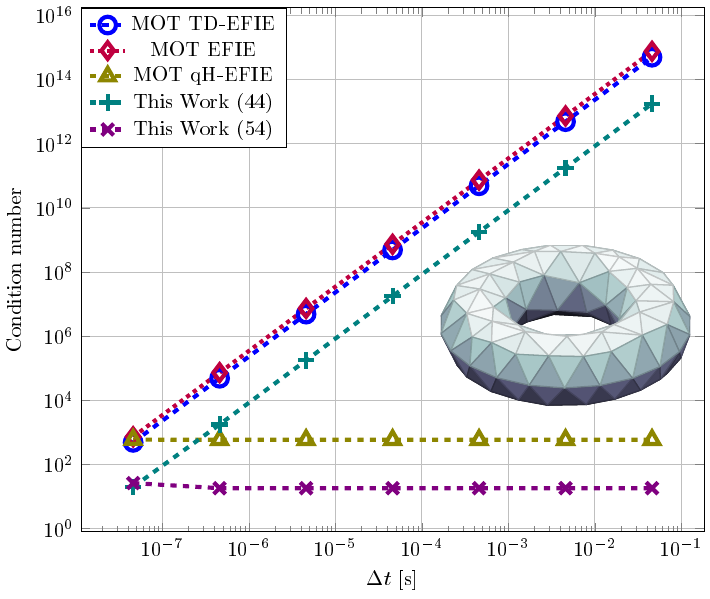}
    \caption{Condition number with respect to the time step size ($N_s=387$) on a torus of the EFIE MOT schemes.}
    \label{fig:torus_cond_number_dt}
\end{figure}

\subsection{Non-canonical geometries}

The final set of numerical tests is dedicated to more complex test structures (\figref{fig::all_geometry}). In addition, instead of direct solver we rely on the iterative solver GMRES with different relative target tolerances $\epsilon$ \cite{saad1986gmres} . For practical reason, the maximum number of iteration has been limited to 200 without restart.
\begin{figure}[ht]
    \centering
    \includegraphics[width=1\linewidth]{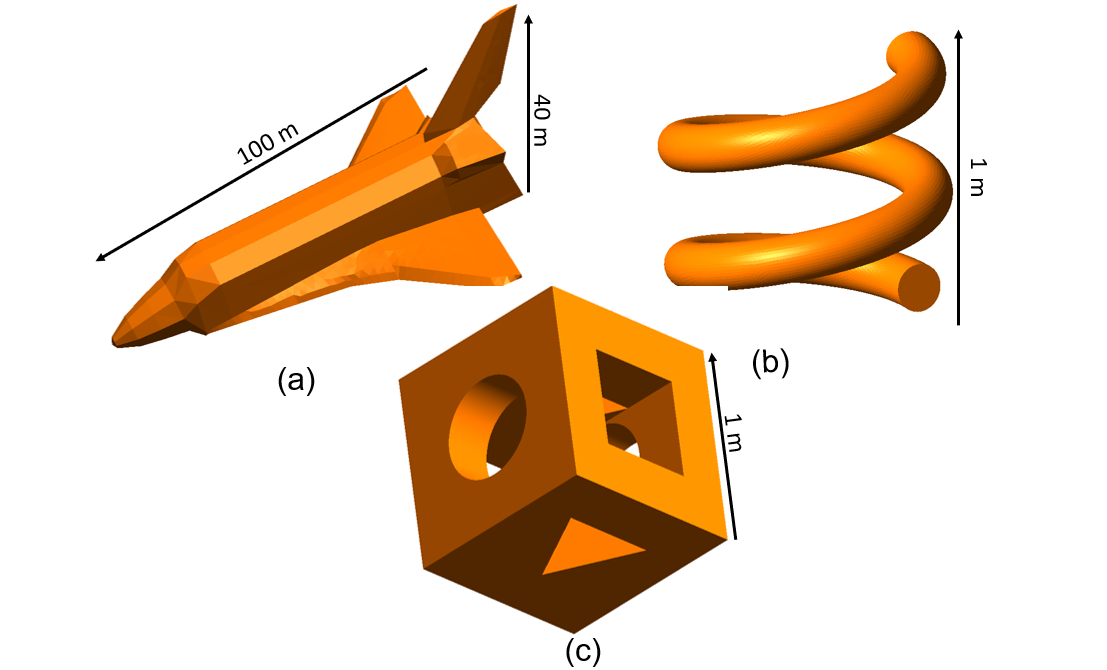}
    \caption{Test structures. (a) space shuttle; (b) compression spring; (c) cube with multiple holes.}
    \label{fig::all_geometry}
\end{figure}

\begin{table}
\begin{center}
\begin{tabular}{|c|c|c|c|c|}
\hline
 \multicolumn{5}{|c|}{ Space shuttle: $N_s=4311$}\\
    \hline
    & Cond($\mat Z_0$)& \multicolumn{3}{|c|}{{\footnotesize $N_\text{iter}$}}\\
    \hline
     Tolerance $\epsilon$ & $\emptyset$ &   $10^{-3}$&    $10^{-6}$ &    $10^{-10}$ \\
    \hline
    MOT TD-EFIE & $4.3\cdot 10^6$ & $>200$ & $>200$ & $>200$  \\    \hline
    MOT EFIE & $6.5\cdot 10^6$ & $>200$ & $>200$ & $>200$  \\    \hline
    MOT qH-EFIE  & $5.3\cdot 10^4$ & $>200$ & $>200$ & $>200$  \\    \hline
    This work \eqref{eq:MOTCPEFIE} & $38$ & $15$ & $28$ & $41$  \\   \hline
    This work \eqref{eq:MOTqHCPEFIE}  & $41$ & $4$ &  $18$ & $33$ \\    \hline
\end{tabular}
\\ \vskip .1cm
\begin{tabular}{|c|c|c|c|c|}
\hline
 \multicolumn{5}{|c|}{Compression spring: $N_s=3906$}\\
    \hline
    & Cond($\mat Z_0$)& \multicolumn{3}{|c|}{{\footnotesize $N_\text{iter}$}}\\
    \hline
     Tolerance $\epsilon$& $\emptyset$ &   $10^{-3}$&    $10^{-6}$ &    $10^{-10}$ \\
    \hline
    MOT TD-EFIE & $2.9\cdot 10^6$ & $>200$ & $>200$ & $>200$  \\    \hline
    MOT EFIE & $4.3\cdot 10^6$ & $>200$ & $>200$ & $>200$  \\    \hline
    MOT qH-EFIE  &$1.6\cdot 10^3$ & $57$ & $154$ & $>200$ \\    \hline
    This work \eqref{eq:MOTCPEFIE} & $158$ & $8$ & $31$ & $50$  \\   \hline
    This work \eqref{eq:MOTqHCPEFIE} & $133$ & $2$ &  $19$ & $39$ \\    \hline
\end{tabular}
\\ \vskip .1cm
\begin{tabular}{|c|c|c|c|c|}
\hline
 \multicolumn{5}{|c|}{Cube with multiple holes: $N_s=3267$}\\
    \hline
    & Cond($\mat Z_0$)& \multicolumn{3}{|c|}{{\footnotesize $N_\text{iter}$}}\\
    \hline
    Tolerance $\epsilon$ & $\emptyset$ &   $10^{-3}$&    $10^{-6}$ &    $10^{-10}$ \\
    \hline
    MOT TD-EFIE & $1.2\cdot 10^9$ & $>200$ & $>200$ & $>200$  \\    \hline
    MOT EFIE & $1.8\cdot 10^9$ & $>200$ & $>200$ & $>200$  \\    \hline
    MOT qH-EFIE  &$6.6\cdot 10^3$ & $50$ & $196$ & $>200$  \\    \hline
    This work \eqref{eq:MOTCPEFIE} & $9.5 \cdot 10^6$ & $5$ & $10$ & $82$  \\   \hline
    This work \eqref{eq:MOTqHCPEFIE} & $50$ & $1$ &  $11$ & $25$ \\    \hline
\end{tabular}
\end{center}
\captionsetup{width=1\textwidth}
\caption{Condition number and maximum number of iterations of the different formulations with different relative tolerances.}
\label{table_shuttle_ring_cube}
\end{table}
All the structures have been illuminated by a pulse Gaussian plane wave with $f_\mr{bw}=1.6\si{\mega \hertz}$. Table~\ref{table_shuttle_ring_cube} shows the condition number and the number of iterations needed. The Calderón preconditioned formulation (CP-EFIE) \eqref{eq:MOTCPEFIE} and the Calderón preconditioned formulation regularized by the quasi-Helmholtz-projectors (qH-CP-EFIE) \eqref{eq:RegCPEFIE} are the only formulations requiring less than 200 iterations to converge at each time steps. As expected, the condition number of the CP-EFIE is high on non-simply connected geometries because of the presence of the operator DC nullspace on these structures. Even if, in this case, the number of iterations remains low, the solution is corrupted by DC instability arising from this nullspace. This phenomenon is absent for the qH-CP-EFIE formulation which is free from high conditioning or DC instability and yields stable solutions up to the target precision of the iterative solver.

\section*{Conclusion}

In this paper, novel Calder\'on preconditioned techniques have been presented for the time domain Electric Field Integral Equations solved with Marching-On-In-Time with convolution quadratures. These formulations eliminate the DC-instability for simply and multiply connected geometries. In addition, they cure the $h$-refinement and large time step breakdowns and generate well-conditioned Marching-On-In Time. Finally, numerical results on complex geometries showcased the effectiveness of the proposed schemes.

\printbibliography

@article{lubich1988convolutionI,
  title={Convolution quadrature and discretized operational calculus. I},
  author={Lubich, Christian},
  journal={Numerische Mathematik},
  volume={52},
  number={2},
  pages={129--145},
  year={1988},
  publisher={Springer}
}

@article{lubich1988convolutionII,
  title={Convolution quadrature and discretized operational calculus. II},
  author={Lubich, Christian},
  journal={Numerische Mathematik},
  volume={52},
  number={4},
  pages={413--425},
  year={1988},
  publisher={Springer}
}

@article{dely2019large,
  title={Large time step and DC stable TD-EFIE discretized with implicit Runge--Kutta methods},
  author={D{\'e}ly, A. and Andriulli, F. P and Cools, K.},
  journal={IEEE Transactions on Antennas and Propagation},
  volume={68},
  number={2},
  pages={976--985},
  year={2019},
  publisher={IEEE}
}

@article{wang2011implicit,
  title={Implicit Runge-Kutta methods for the discretization of time domain integral equations},
  author={Wang, X. and Weile, D. S},
  journal={IEEE transactions on antennas and propagation},
  volume={59},
  number={12},
  pages={4651--4663},
  year={2011},
  publisher={IEEE}
}

@ARTICLE{9580445,
  author={Adrian, Simon B. and Dély, Alexandre and Consoli, Davide and Merlini, Adrien and Andriulli, Francesco P.},
  journal={IEEE Open Journal of Antennas and Propagation}, 
  title={Electromagnetic Integral Equations: Insights in Conditioning and Preconditioning}, 
  year={2021},
  volume={2},
  number={},
  pages={1143-1174},
  doi={10.1109/OJAP.2021.3121097}}

@article{banjai2022integral,
  title={Integral equation methods for evolutionary PDE},
  author={Banjai, L and Sayas, FJ},
  journal={Springer Ser. Comput. Math., Springer, Cham},
  year={2022}
}

@article{adrian2021electromagnetic,
  title={Electromagnetic Integral Equations: Insights in Conditioning and Preconditioning},
  author={Adrian, Simon B and Dély, Alexandre and Consoli, Davide and Merlini, Adrien and Andriulli, Francesco P},
  journal={IEEE Open Journal of Antennas and Propagation},
  year={2021},
  publisher={IEEE}
}

@article{cools2009nullspaces,
  title={Nullspaces of MFIE and Calder{\'o}n preconditioned EFIE operators applied to toroidal surfaces},
  author={Cools, Kristof and Andriulli, Francesco P and Olyslager, Femke and Michielssen, Eric},
  journal={IEEE Transactions on Antennas and Propagation},
  volume={57},
  number={10},
  pages={3205--3215},
  year={2009},
  publisher={IEEE}
}

@article{buffa2007dual,
  title={A dual finite element complex on the barycentric refinement},
  author={Buffa, Annalisa and Christiansen, Snorre},
  journal={Mathematics of Computation},
  volume={76},
  number={260},
  pages={1743--1769},
  year={2007}
}

@book{bennett1968technique,
  title={A technique for computing approximate electromagnetic impulse response of conducting bodies},
  author={Bennett Jr, C Leonard},
  year={1968},
  publisher={Purdue University}
}

@article{lubich1993runge,
  title={Runge-Kutta methods for parabolic equations and convolution quadrature},
  author={Lubich, Ch and Ostermann, A},
  journal={mathematics of computation},
  volume={60},
  number={201},
  pages={105--131},
  year={1993}
}

@book{oppenheim2001discrete,
  title={Discrete-time signal processing. Vol. 2},
  author={Oppenheim, Alan V and Buck, John R and Schafer, Ronald W},
  year={2001},
  publisher={Upper Saddle River, NJ: Prentice Hall}
}

@article{schanz1997new,
  title={A new visco-and elastodynamic time domain boundary element formulation},
  author={Schanz, Martin and Antes, Heinz},
  journal={Computational Mechanics},
  volume={20},
  number={5},
  pages={452--459},
  year={1997},
  publisher={Springer}
}

@article{banjai2012runge,
  title={Runge--Kutta convolution quadrature for the boundary element method},
  author={Banjai, Lehel and Messner, Matthias and Schanz, Martin},
  journal={Computer methods in applied mechanics and engineering},
  volume={245},
  pages={90--101},
  year={2012},
  publisher={Elsevier}
}

@book{schwartz2008mathematics,
  title={Mathematics for the physical sciences},
  author={Schwartz, Laurent},
  year={2008},
  publisher={Courier Dover Publications}
}

@book{butcher2016numerical,
  title={Numerical methods for ordinary differential equations},
  author={Butcher, John Charles},
  year={2016},
  publisher={John Wiley \& Sons}
}

@article{skvortsov2003accuracy,
  title={Accuracy of Runge-Kutta methods applied to stiff problems},
  author={Skvortsov, LM},
  journal={Computational mathematics and mathematical physics},
  volume={43},
  number={9},
  pages={1320--1330},
  year={2003},
  publisher={Oxford; New York: Pergamon Press, c1992-}
}

@article{runge1895numerische,
  title={{\"U}ber die numerische Aufl{\"o}sung von Differentialgleichungen},
  author={Runge, Carl},
  journal={Mathematische Annalen},
  volume={46},
  number={2},
  pages={167--178},
  year={1895},
  publisher={Springer}
}

@article{kutta1901beitrag,
  title={Beitrag zur naherungsweisen integration totaler differentialgleichungen},
  author={Kutta, Wilhelm},
  journal={Z. Math. Phys.},
  volume={46},
  pages={435--453},
  year={1901}
}

@article{rao1982electromagnetic,
  title={Electromagnetic scattering by surfaces of arbitrary shape},
  author={Rao, Sadasiva and Wilton, Donald and Glisson, Allen},
  journal={IEEE Transactions on antennas and propagation},
  volume={30},
  number={3},
  pages={409--418},
  year={1982},
  publisher={IEEE}
}

@incollection{raviart1977mixed,
  title={A mixed finite element method for 2-nd order elliptic problems},
  author={Raviart, Pierre-Arnaud and Thomas, Jean-Marie},
  booktitle={Mathematical aspects of finite element methods},
  pages={292--315},
  year={1977},
  publisher={Springer}
}

@article{cools2011accurate,
  title={Accurate and conforming mixed discretization of the MFIE},
  author={Cools, Kristof and Andriulli, FP and De Zutter, Dani{\"e}l and Michielssen, Eric},
  journal={IEEE antennas and wireless propagation letters},
  volume={10},
  pages={528--531},
  year={2011},
  publisher={IEEE}
}

@article{maruyama2016transient,
  title={Transient elastic wave analysis of 3-D large-scale cavities by fast multipole BEM using implicit Runge--Kutta convolution quadrature},
  author={Maruyama, T and Saitoh, T and Bui, TQ and Hirose, S},
  journal={Computer Methods in Applied Mechanics and Engineering},
  volume={303},
  pages={231--259},
  year={2016},
  publisher={Elsevier}
}

@article{banjai2014fast,
  title={Fast convolution quadrature for the wave equation in three dimensions},
  author={Banjai, Lehel and Kachanovska, Maryna},
  journal={Journal of Computational Physics},
  volume={279},
  pages={103--126},
  year={2014},
  publisher={Elsevier}
}

@book{colton2013integral,
  title={Integral equation methods in scattering theory},
  author={Colton, David and Kress, Rainer},
  year={2013},
  publisher={SIAM}
}

@article{walker2002stability,
  title={The stability of integral equation time-domain scattering computations for three-dimensional scattering; similarities and differences between electrodynamic and elastodynamic computations},
  author={Walker, SP and Bluck, MJ and Chatzis, I},
  journal={International Journal of Numerical Modelling: Electronic Networks, Devices and Fields},
  volume={15},
  number={5-6},
  pages={459--474},
  year={2002},
  publisher={Wiley Online Library}
}

@book{peterson1998computational,
  title={Computational methods for electromagnetics},
  author={Peterson, Andrew F and Ray, Scott L and Mittra, Raj and Institute of Electrical and Electronics Engineers},
  volume={351},
  year={1998},
  publisher={IEEE press New York}
}

@book{jin2011theory,
  title={Theory and computation of electromagnetic fields},
  author={Jin, Jian-Ming},
  year={2011},
  publisher={John Wiley \& Sons}
}

@book{gibson2021method,
  title={The method of moments in electromagnetics},
  author={Gibson, Walton C},
  year={2021},
  publisher={Chapman and Hall/CRC}
}

@book{golub2013matrix,
  title={Matrix computations},
  author={Golub, Gene H and Van Loan, Charles F},
  year={2013},
  publisher={JHU press}
}

@article{rao1991transient,
  title={Transient scattering by conducting surfaces of arbitrary shape},
  author={Rao, Sadasiva M and Wilton, Donald R},
  journal={IEEE Transactions on Antennas and Propagation},
  volume={39},
  number={1},
  pages={56--61},
  year={1991},
  publisher={IEEE}
}

@article{shi2014static,
  title={On the static loop modes in the marching-on-in-time solution of the time-domain electric field integral equation},
  author={Shi, Yifei and Ba{\u{g}}c{\i}, Hakan and Lu, Mingyu},
  journal={IEEE Antennas and Wireless Propagation Letters},
  volume={13},
  pages={317--320},
  year={2014},
  publisher={IEEE}
}

@article{andriulli2009analysis,
  title={Analysis and regularization of the TD-EFIE low-frequency breakdown},
  author={Andriulli, Francesco P and Bagci, Hakan and Vipiana, Francesca and Vecchi, Giuseppe and Michielssen, Eric},
  journal={IEEE transactions on antennas and propagation},
  volume={57},
  number={7},
  pages={2034--2046},
  year={2009},
  publisher={IEEE}
}

@article{chen2001integral,
  title={Integral-equation-based analysis of transient scattering and radiation from conducting bodies at very low frequencies},
  author={Chen, N-W and Ayg{\"u}n, K and Michielssen, E},
  journal={IEE Proceedings-Microwaves, Antennas and Propagation},
  volume={148},
  number={6},
  pages={381--387},
  year={2001},
  publisher={IET}
}

@article{andriulli2012well,
  title={On a well-conditioned electric field integral operator for multiply connected geometries},
  author={Andriulli, Francesco P and Cools, Kristof and Bogaert, Ignace and Michielssen, Eric},
  journal={IEEE transactions on antennas and propagation},
  volume={61},
  number={4},
  pages={2077--2087},
  year={2012},
  publisher={IEEE}
}

@article{bogaert2013low,
  title={Low-frequency scaling of the standard and mixed magnetic field and M{\"u}ller integral equations},
  author={Bogaert, Ignace and Cools, Kristof and Andriulli, Francesco P and Ba{\u{g}}c{\i}, Hakan},
  journal={IEEE Transactions on Antennas and Propagation},
  volume={62},
  number={2},
  pages={822--831},
  year={2013},
  publisher={IEEE}
}

@article{andriulli2009time,
  title={Time domain Calder{\'o}n identities and their application to the integral equation analysis of scattering by PEC objects part II: Stability},
  author={Andriulli, Francesco P and Cools, Kristof and Olyslager, Femke and Michielssen, Eric},
  journal={IEEE transactions on antennas and propagation},
  volume={57},
  number={8},
  pages={2365--2375},
  year={2009},
  publisher={IEEE}
}

@article{cools2009time,
  title={Time domain Calder{\'o}n identities and their application to the integral equation analysis of scattering by PEC objects Part I: Preconditioning},
  author={Cools, Kristof and Andriulli, Francesco P and Olyslager, Femke and Michielssen, Eric},
  journal={IEEE Transactions on Antennas and Propagation},
  volume={57},
  number={8},
  pages={2352--2364},
  year={2009},
  publisher={IEEE}
}

@article{andriulli2008multiplicative,
  title={A multiplicative Calderon preconditioner for the electric field integral equation},
  author={Andriulli, Francesco P and Cools, Kristof and Bagci, Hakan and Olyslager, Femke and Buffa, Annalisa and Christiansen, Snorre and Michielssen, Eric},
  journal={IEEE Transactions on Antennas and Propagation},
  volume={56},
  number={8},
  pages={2398--2412},
  year={2008},
  publisher={IEEE}
}

@article{beghein2015dc,
  title={A DC-stable, well-balanced, Calder{\'o}n preconditioned time domain electric field integral equation},
  author={Beghein, Yves and Cools, Kristof and Andriulli, Francesco P},
  journal={IEEE Transactions on Antennas and Propagation},
  volume={63},
  number={12},
  pages={5650--5660},
  year={2015},
  publisher={IEEE}
}

@article{chipman1971stable,
  title={A-stable Runge-Kutta processes},
  author={Chipman, FH},
  journal={BIT Numerical Mathematics},
  volume={11},
  number={4},
  pages={384--388},
  year={1971},
  publisher={Springer}
}

@inproceedings{andriulli2008dottrick,
  title={The “dottrick TDEFIE”: a DC stable integral equation for analyzing transient scattering from PEC bodies},
  author={Andriulli, Francesco P and Cools, Kristof and Olyslager, Femke and Michielssen, Eric},
  booktitle={2008 IEEE Antennas and Propagation Society International Symposium},
  pages={1--4},
  year={2008},
  organization={IEEE}
}

@inproceedings{andriulli2012analysis,
  title={Analysis and discretization of the Yukawa-Calderon preconditioned CFIE},
  author={Andriulli, Francesco P and Cools, Kristof and Bogaert, Ignace and Bagci, Hakan and Yla-Ojiala, Pasi and Michielssen, Eric},
  booktitle={28th Annual Review of Progress in Applied Computational Electromagnetics (ACES-2012)},
  pages={454--463},
  year={2012},
  organization={Curran Associates, Inc.}
}

@book{harrington1996field,
  title={Field computation by moment methods},
  author={Harrington, Roger F and Harrington, Jan L},
  year={1996},
  publisher={Oxford University Press, Inc.}
}

@book{sayas2016retarded,
  title={Retarded potentials and time domain boundary integral equations: A road map},
  author={Sayas, Francisco-Javier},
  volume={50},
  year={2016},
  publisher={Springer}
}

@book{nedelec2001acoustic,
  title={Acoustic and electromagnetic equations: integral representations for harmonic problems},
  author={N{\'e}d{\'e}lec, Jean-Claude},
  volume={144},
  year={2001},
  publisher={Springer}
}

@article{chung2004solution,
  title={Solution of time domain electric field integral equation using the Laguerre polynomials},
  author={Chung, Young-Seek and Sarkar, Tapan K and Jung, Baek Ho and Salazar-Palma, Magdalena and Ji, Zhong and Jang, Seongman and Kim, Kyungjung},
  journal={IEEE Transactions on Antennas and Propagation},
  volume={52},
  number={9},
  pages={2319--2328},
  year={2004},
  publisher={IEEE}
}

@inproceedings{beghein2013temporal,
  title={A temporal Galerkin discretization of the charge-current continuity equation},
  author={Beghein, Yves and Cools, Kristof and De Zutter, Dani{\"e}l},
  booktitle={2013 International Conference on Electromagnetics in Advanced Applications (ICEAA)},
  pages={628--631},
  year={2013},
  organization={IEEE}
}

@article{hsiao1997mathematical,
  title={Mathematical foundations for error estimation in numerical solutions of integral equations in electromagnetics},
  author={Hsiao, George C and Kleinman, Ralph E},
  journal={IEEE transactions on Antennas and Propagation},
  volume={45},
  number={3},
  pages={316--328},
  year={1997},
  publisher={IEEE}
}

@article{van2022role,
  title={The Role of Jordan Blocks in the MOT-Scheme Time Domain EFIE Linear-in-Time Solution Instability},
  author={van Diepen, Petrus WN and Dilz, Roeland J and Zwamborn, Adrianus PM and van Beurden, Martijn C},
  journal={Progress In Electromagnetics Research B},
  volume={95},
  number={95},
  pages={123},
  year={2022},
  publisher={Electromagnetics Academy}
}

@inproceedings{cordel2022calderon,
  title={Calder{\'o}n Preconditioners for the TD-EFIE discretized with Convolution Quadratures},
  author={Cordel, Pierrick and D{\'e}ly, Alexandre and Merlini, Adrien and Andriulli, Francesco P},
  booktitle={2022 IEEE International Symposium on Antennas and Propagation and USNC-URSI Radio Science Meeting (AP-S/URSI)},
  pages={1672--1673},
  year={2022},
  organization={IEEE}
}

@article{saad1986gmres,
  title={GMRES: A generalized minimal residual algorithm for solving nonsymmetric linear systems},
  author={Saad, Youcef and Schultz, Martin H},
  journal={SIAM Journal on scientific and statistical computing},
  volume={7},
  number={3},
  pages={856--869},
  year={1986},
  publisher={SIAM}
}

@article{dely2022convolution,
  title={Convolution Quadrature Time Domain Integral Equation Methods for Electromagnetic Scattering},
  author={D{\'e}ly, Alexandre and Merlini, Adrien and Cools, Kristof and Andriulli, Francesco P},
  journal={Advances in Time-Domain Computational Electromagnetic Methods},
  pages={321--359},
  year={2022},
  publisher={Wiley Online Library}
}

\begin{IEEEbiography}[{\includegraphics[width=1in,height=1.25in,clip,keepaspectratio]{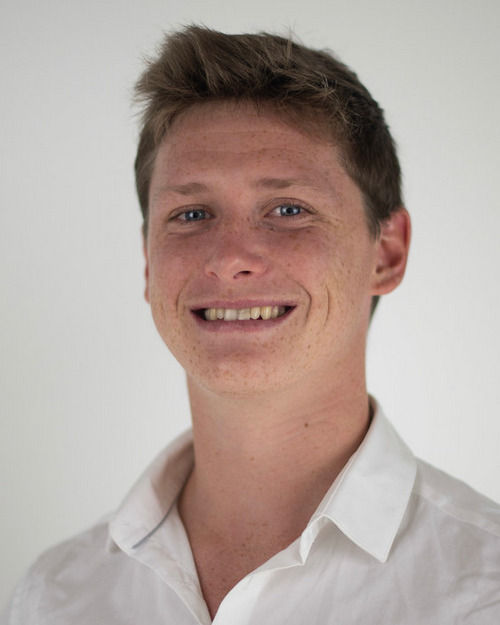}}]{PIERRICK CORDEL } (Student Member, IEEE) received the M.Sc. Eng. Degree from the Ecole Nationale Supérieure des Mines de Nancy, France, in 2021 as well as a M.Sc. in Industrial Mathematics from the University of Luxembourg, Luxembourg, the same year. Currently, he is doing a PhD thesis at the institute Politecnico di Torino, in Italy. 

His research focuses on time domain integral equations discretized with the convolution quadrature method.
\end{IEEEbiography}
\begin{IEEEbiography}[{\includegraphics[width=1in,height=1.25in,clip,keepaspectratio]{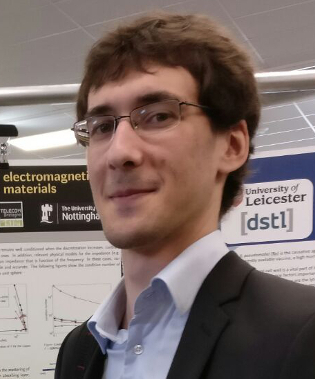}}]{ALEXANDRE DÉLY } 
received the M.Sc. Eng. degree from the École Nationale Supérieure des Télécommunications de Bretagne (Télécom Bretagne), France, in 2015. He received the Ph.D. degree from the École Nationale Supérieure Mines-Télécom Atlantique (IMT Atlantique), France, and from the University of Nottingham, United Kingdom, in 2019.

His research focuses on preconditioned and fast solution of boundary element methods frequency domain and time domain integral equations. He is currently working in Thales, Elancourt, France, on electromagnetic modeling and numerical simulations
\end{IEEEbiography}

\begin{IEEEbiography}[{\includegraphics[width=1in,height=1.25in,clip,keepaspectratio]{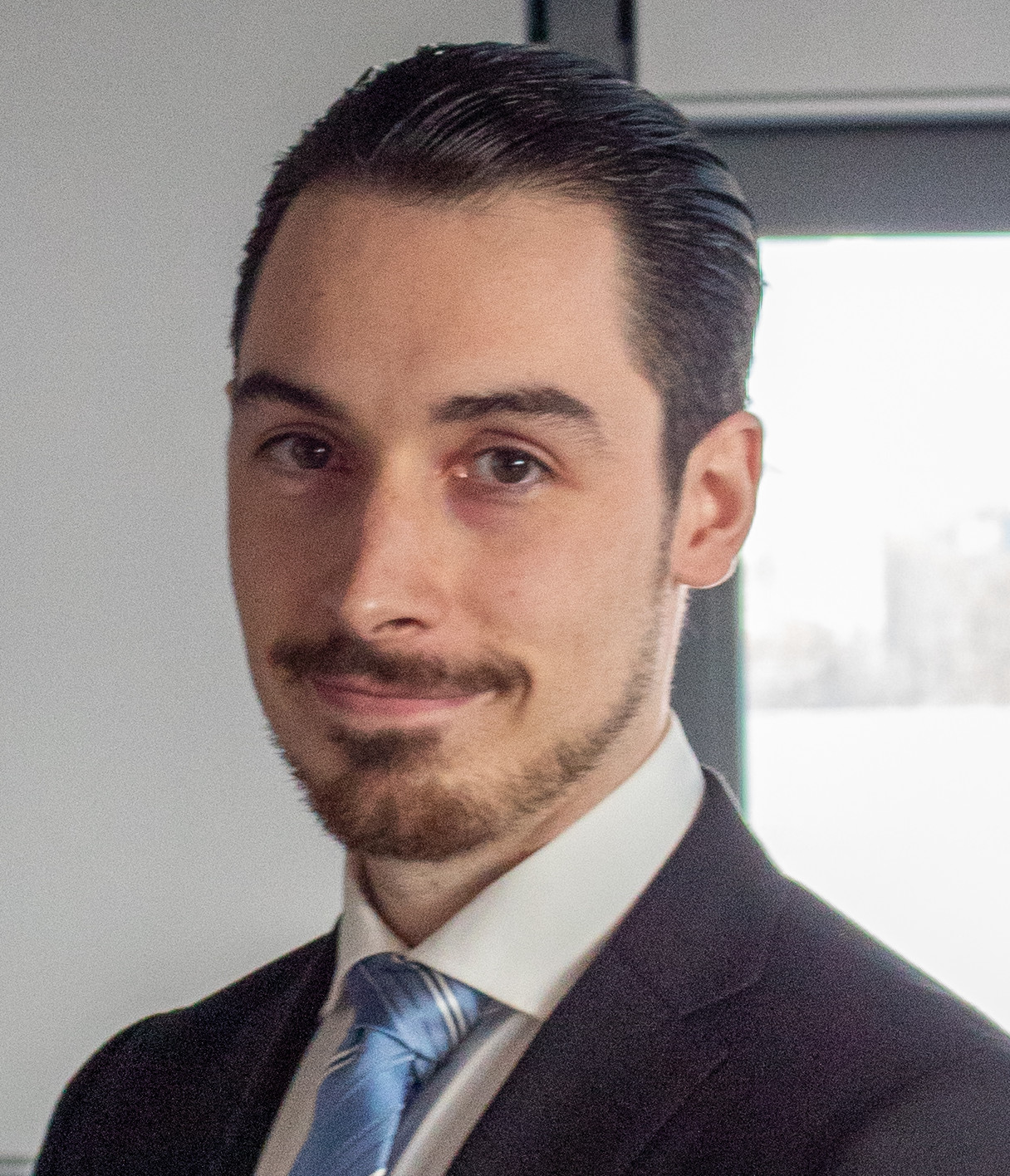}}]{ADRIEN MERLINI } 
(S’16–M’19) received the M.Sc.Eng. degree from the École Nationale Supérieure des Télécommunications de Bretagne (Télécom Bretagne), France, in 2015 and received the Ph.D. degree from the École Nationale Supérieure MinesTélécom Atlantique (IMT Atlantique), France, in 2019.
From 2018 to 2019, he was a visiting Ph.D. student at the Politecnico di Torino, Italy, which he then joined as a Research Associate. Since 2019, he has been an Associate Professor with the Microwave Department, IMT Atlantique. His research interests include preconditioning and acceleration of integral equation solvers for electromagnetic simulations and their application in brain imaging. 

Dr. Merlini received 2 Young Scientist Awards at the URSI GASS 2020 and the EMTS 2023 meetings. In addition, he has co-authored a paper that received the 2022 ICEAA-IEEE APWC Best Paper Award, 5 that received honorable mentions (URSI/IEEE–APS 2021, 2022, and 2023) and 3 best paper finalists (URSI GASS 2020, URSI/IEEE–APS 2021 and 2022). He is a member of IEEE-HKN, the IEEE Antennas and Propagation Society, URSI France, and of the Lab-STICC laboratory. He is currently serving as Associate Editor for the Antenna and Propagation Magazine.
\end{IEEEbiography}

\begin{IEEEbiography}[{\includegraphics[width=1in,height=1.25in,clip,keepaspectratio]{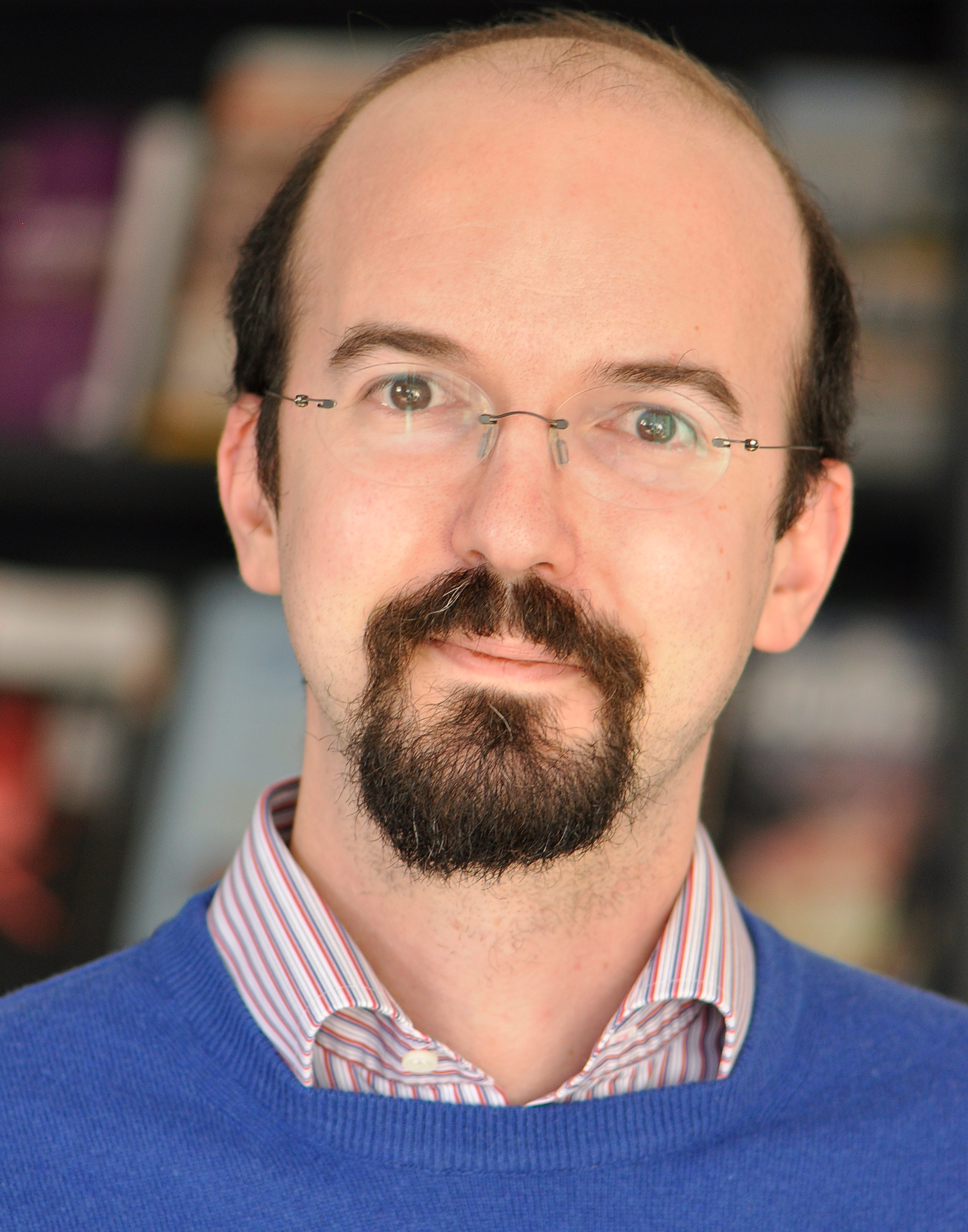}}]{FRANCESCO P. ANDRIULLI } (S’05–M’09– SM’11–F’23) received the Laurea in electrical engineering from the Politecnico di Torino, Italy, in 2004, the MSc in electrical engineering and computer science from the University of Illinois at Chicago in 2004, and the PhD in electrical engineering from the University of Michigan at Ann Arbor in 2008. From 2008 to 2010 he was a Research Associate with the Politecnico di Torino. From 2010 to 2017 he was an Associate Professor (2010-2014) and then Full Professor with the École Nationale Supérieure Mines-Télécom Atlantique (IMT Atlantique, previously ENST Bretagne), Brest, France. Since 2017 he has been a Full Professor with the Politecnico di Torino, Turin, Italy. His research interests are in computational electromagnetics with focus on frequency- and time-domain integral equation solvers, well-conditioned formulations, fast solvers, low-frequency electromagnetic analyses, and modeling techniques for antennas, wireless components, microwave circuits, and biomedical applications with a special focus on brain imaging.

Prof. Andriulli received several best paper awards at conferences and symposia (URSI NA 2007, IEEE AP-S 2008, ICEAA IEEE-APWC 2015) also in co-authorship with his students and collaborators (ICEAA IEEE-APWC 2021, EMTS 2016, URSI-DE Meeting 2014, ICEAA 2009) with whom received also a second prize conference paper (URSI GASS 2014), a third prize conference paper (IEEE–APS 2018), seven honorable mention conference papers (ICEAA 2011, URSI/IEEE–APS 2013, 4 in URSI/IEEE–APS 2022, URSI/IEEE–APS 2023) and other three finalist conference papers (URSI/IEEE-APS 2012, URSI/IEEE-APS 2007, URSI/IEEE-APS 2006, URSI/IEEE–APS 2022). A Fellow of the IEEE, he is also the recipient of the 2014 IEEE AP-S Donald G. Dudley Jr. Undergraduate Teaching Award, of the triennium 2014-2016 URSI Issac Koga Gold Medal, and of the 2015 L. B. Felsen Award for Excellence in Electrodynamics. 

Prof. Andriulli is a member of Eta Kappa Nu, Tau Beta Pi, Phi Kappa Phi, and of the International Union of Radio Science (URSI). He is the Editorin-Chief of the IEEE Antennas and Propagation Magazine and he serves as a Track Editor for the IEEE Transactions on Antennas and Propagation, and as an Associate Editor for URSI Radio Science Letters. He served as an Associate Editor for the IEEE Transactions on Antennas and Propagation, IEEE Antennas and Wireless Propagation Letters, IEEE Access and IET-MAP.
\end{IEEEbiography}

\end{document}